\documentclass[12pt]{article}
\usepackage[letterpaper, left=2.5cm, right=2.5cm, top=2.5cm,bottom=2.5cm,dvips]{geometry}\usepackage{epsf,amssymb,latexsym,amsmath}
\usepackage{graphicx}
\usepackage{tikz}
\usetikzlibrary{decorations.pathreplacing}
\usepackage[autostyle=true]{csquotes}

\newcommand{\proof}{\noindent{\bf Proof.\ }}
\newcommand{\qed}{\hfill $\square$ \bigskip}

\newtheorem{theorem}{\bf Theorem}[section]
\newtheorem{corollary}[theorem]{\bf Corollary}
\newtheorem{lemma}[theorem]{\bf Lemma}
\newtheorem{proposition}[theorem]{\bf Proposition}
\newtheorem{conjecture}[theorem]{\bf Conjecture}

\newtheorem{definition}[theorem]{\bf Definition}
\newtheorem{question}[theorem]{\bf Question}
\newtheorem{problem}[theorem]{\bf Problem}

\newcommand{\ecc}{{\rm ecc}}

\begin{document}

\title{New results on   Wiener index of trees with a given diameter}

\author{
Bojana Borovi\' canin \\
\small \it Faculty of Science\\
\small \it University of Kragujevac, Serbia \\
\small \it e-mail: bojana.borovicanin@pmf.kg.ac.rs \\
\and
Dragana Bo\v zovi\'{c} \\
\small \it Faculty of Electrical Engineering \\
\small \it and Computer Science \\
\small \it University of Maribor, Slovenia \\
\small \it e-mail: dragana.bozovic@um.si\\
\and
Edin Glogi\' c \\ 
\small \it State University of Novi Pazar \\
\small \it e-mail: edinglogic@np.ac.rs\\
\and
Daša Mesari\v c \v Stesl \\
\small \it Faculty of Computer and Information Science \\
\small \it University of Ljubljana, Slovenia \\
\small \it e-mail: dasa.stesl@fri.uni-lj.si\\
\and
Simon \v Spacapan \\
\small \it Faculty of Mechanical Engineering \\
\small \it University of Maribor, Slovenia \\
\small \it e-mail: simon.spacapan@um.si\\
\and
Emir Zogi\' c \\
\small \it State University of Novi Pazar \\
\small \it e-mail: ezogic@np.ac.rs
}
\date{}

\maketitle

\begin{abstract} 
\noindent We study the Wiener index of a class of trees with  fixed diameter and order. 
A double broom is a tree such that there exist two vertices $u$ and $v$, such that each leaf of $T$ is adjacent to $u$ or $v$. 
We prove that for a tree $T$ of diameter $d$ and (sufficiently large) order $n$ such that $n\leq d-2+4 \left\lfloor \sqrt{ \frac{d-1}{2}} \right\rfloor$, 
$T$ has maximum Wiener index (in the class of trees of diameter $d$ and order $n$) if and only if $T$ is a balanced double broom. Our results are sharp up to a small constant.  
\end{abstract}



\section{Introduction}


 Wiener index is a topological index in graph theory which was introduced in 1947 by Harry Wiener in  \cite{wiener}. It has been extensively studied by various authors (see e.g. \cite{DoEnGu, DoGuKlZi, DoMe, knor, Xu}) because of its applicability in network analysis, combinatorics, and chemistry. 

All graphs considered in this article are finite, simple and undirected. For a graph $G$, we denote its vertex set by $V(G)$ and its edge set by $E(G)$. 
The distance $d(u, v)$ between   vertices $u, v \in V(G)$ is the length of a shortest $uv$-path in $G$.
A graph is connected if for any pair of its vertices there exists a path between them. The diameter of a graph $G$, denoted as $diam(G)$, is the 
maximum distance between a pair of vertices of $G$. For a connected graph $G$, the \textit{eccentricity} of $v\in V(G)$, denoted as $\ecc(v)$, is the maximum distance between $v$ and a vertex of $G$. 
For a positive integer $n$ we denote $[n]=\{1,\ldots, n\}$ and $[n]_0=[n]\cup\{0\}$.

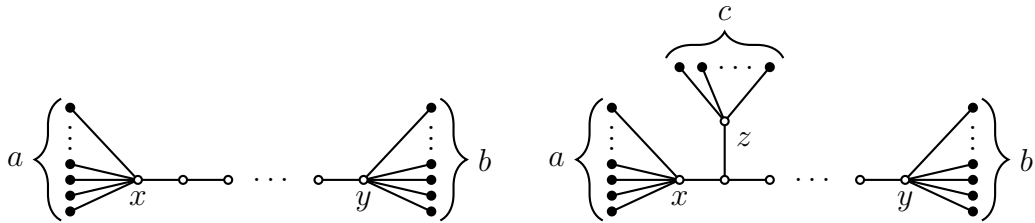
\begin{figure}[hbt!]
	\label{broom}
	\begin{center}
		\begin{tikzpicture}[scale=.6,style=thick,x=1cm,y=1cm]
			\def\vr{2.5pt} 
			\path (-0.5,0) coordinate (v1);
			\path (1,0) coordinate (v2);
			\path (2,0) coordinate (v3);
			\path (3,0) coordinate (v4);
			\path (5,0) coordinate (vd2);
			\path (6,0) coordinate (vd1);
			\path (7.5,0) coordinate (vd);
			\node at (4,0){$\dots$};
			
			\path (-0.5,-0.7) coordinate (v5);
			\path (-0.5,-0.35) coordinate (v6);
			\path (-0.5,0.35) coordinate (v7);
			\node at (-0.5,1.1){$\vdots$};
			\path (-0.5,1.6) coordinate (v8);
			
			\draw[decorate,decoration={brace,amplitude=10pt,mirror}]
			(-0.7,1.8) -- (-0.7,-0.9) node[midway, left=10pt] {$a$};
			
			\path (7.5,-0.7) coordinate (u5);
			\path (7.5,-0.35) coordinate (u6);
			\path (7.5,0.35) coordinate (u7);
			\node at (7.5,1.1){$\vdots$};
			\path (7.5,1.6) coordinate (u8);
			
			\draw[decorate,decoration={brace,amplitude=10pt}]
			(7.7,1.8) -- (7.7,-0.9) node[midway, right=10pt] {$b$};
			
			
			\draw (v1) -- (v2);
			\draw (v2) -- (v3);
			\draw (v3) -- (v4);
			\draw (vd2) -- (vd1);
			\draw (vd1) -- (vd);
			\draw (v2) -- (v5);
			\draw (v2) -- (v6);
			\draw (v2) -- (v7);
			\draw (v2) -- (v8);
			\draw (vd1) -- (u5);
			\draw (vd1) -- (u6);
			\draw (vd1) -- (u7);
			\draw (vd1) -- (u8);
			
			\draw (v1) [fill=black] circle (\vr);
			\draw (v2) [fill=white] circle (\vr);
			\draw (v3) [fill=white] circle (\vr);
			\draw (v4) [fill=white] circle (\vr);
			\draw (vd2) [fill=white] circle (\vr);
			\draw (vd1) [fill=white] circle (\vr);
			\draw (vd) [fill=black] circle (\vr);
			\draw (v5) [fill=black] circle (\vr);
			\draw (v6) [fill=black] circle (\vr);
			\draw (v7) [fill=black] circle (\vr);
			\draw (v8) [fill=black] circle (\vr);
			\draw (u5) [fill=black] circle (\vr);
			\draw (u6) [fill=black] circle (\vr);
			\draw (u7) [fill=black] circle (\vr);
			\draw (u8) [fill=black] circle (\vr);
			
			\draw[anchor = north] (v2) node {$x$};
			\draw[anchor = north] (vd1) node {$y$};
			
			\path (11.5,0) coordinate (z1);
			\path (13,0) coordinate (z2);
			\path (14,0) coordinate (z3);
			\path (15,0) coordinate (z4);
			\path (17,0) coordinate (zd2);
			\path (18,0) coordinate (zd1);
			\path (19.5,0) coordinate (zd);
			\node at (16,0){$\dots$};
			
			\path (14,1.3) coordinate (z);
			
			\path (11.5,-0.7) coordinate (z5);
			\path (11.5,-0.35) coordinate (z6);
			\path (11.5,0.35) coordinate (z7);
			\node at (11.5,1.1){$\vdots$};
			\path (11.5,1.6) coordinate (z8);
			
			\draw[decorate,decoration={brace,amplitude=10pt,mirror}]
			(11.3,1.8) -- (11.3,-0.9) node[midway, left=10pt] {$a$};
			
			\path (19.5,-0.7) coordinate (w5);
			\path (19.5,-0.35) coordinate (w6);
			\path (19.5,0.35) coordinate (w7);
			\node at (19.5,1.1){$\vdots$};
			\path (19.5,1.6) coordinate (w8);
			
			\draw[decorate,decoration={brace,amplitude=10pt}]
			(19.7,1.8) -- (19.7,-0.9) node[midway, right=10pt] {$b$};
			
			\path (13,2.5) coordinate (g2);
			\path (13.5,2.5) coordinate (g3);
			\node at (14.3,2.5){$\dots$};
			\path (15,2.5) coordinate (g);
			
			\draw[decorate,decoration={brace,amplitude=10pt}]
			(12.7,2.7) -- (15.3,2.7) node[midway, above=10pt] {$c$};
			
			
			\draw (z1) -- (z2);
			\draw (z2) -- (z3);
			\draw (z3) -- (z4);
			\draw (zd2) -- (zd1);
			\draw (zd1) -- (zd);
			\draw (z2) -- (z5);
			\draw (z2) -- (z6);
			\draw (z2) -- (z7);
			\draw (z2) -- (z8);
			\draw (zd1) -- (w5);
			\draw (zd1) -- (w6);
			\draw (zd1) -- (w7);
			\draw (zd1) -- (w8);
			\draw (z) -- (z3);
			\draw (z) -- (g2);
			\draw (z) -- (g3);
			\draw (z) -- (g);
			
			\draw (z1) [fill=black] circle (\vr);
			\draw (z2) [fill=white] circle (\vr);
			\draw (z3) [fill=white] circle (\vr);
			\draw (z4) [fill=white] circle (\vr);
			\draw (zd2) [fill=white] circle (\vr);
			\draw (zd1) [fill=white] circle (\vr);
			\draw (zd) [fill=black] circle (\vr);
			\draw (z5) [fill=black] circle (\vr);
			\draw (z6) [fill=black] circle (\vr);
			\draw (z7) [fill=black] circle (\vr);
			\draw (z8) [fill=black] circle (\vr);
			\draw (w5) [fill=black] circle (\vr);
			\draw (w6) [fill=black] circle (\vr);
			\draw (w7) [fill=black] circle (\vr);
			\draw (w8) [fill=black] circle (\vr);
			
			\draw (z) [fill=white] circle (\vr);
			\draw (g2) [fill=black] circle (\vr);
			\draw (g3) [fill=black] circle (\vr);
			\draw (g) [fill=black] circle (\vr);
			
			\draw[anchor = north] (z2) node {$x$};
			\draw[anchor = north] (zd1) node {$y$};
			\draw[anchor = north west] (z) node {$z$};
			
		\end{tikzpicture}
	\end{center}
	\caption{Double broom (left) and triple broom (right).}
\end{figure}

A {\em tree} is a connected graph without cycles. A vertex of degree one in a tree is called a \textit{leaf}, and a \textit{broom vertex} in a tree $T$ is any vertex adjacent to a leaf of $T$.  A {\em double broom} $B(n, a, b)$ is a tree on $n$ vertices with exactly two broom vertices $x$ and $y$ such that 
$\deg(x)=a+1$ and $\deg(y)=b+1$ and a {\em triple broom}    $B(n,a,b,c)$ is a tree on $n$ vertices with exactly three broom vertices 
$x,y$ and $z$ such that $a,b$ and $c$ leaves are attached to  them, respectively, and such that $d(x,z)=2, d(x,y)=d(z,y)=d-2$, where 
$d=n-a-b-c$ is the diameter of $B(n,a,b,c)$, see Fig.~\ref{broom}. A double broom $B(n,a,b)$ is {\em balanced} if $|a-b|\leq 1$.

For a connected graph $G$, the Wiener index of $G$ is denoted as $W(G)$, and is the sum of distances between all unordered pairs of vertices
$$W(G) = \sum_{\{u,v\} \subseteq V(G)} d(u,v).$$

A central  objective in the study of Wiener index is to determine graphs (in a given class of graphs) that maximize or minimize Wiener index. 
In this respect the classes of graphs (with fixed order) with   fixed radius  \cite{cambie2,ChenWuAn,DasNad,StMiVu,YouLiu}, minimum degree \cite{Alo}, maximum degree \cite{Alo, BozVuk, Dong, Fish, stevanovic}, or diameter have been considered.




Since our work focuses on trees with fixed order and diameter, we now give a more detailed overview of related results in this setting. 
In \cite{Plesnik} Plesn\'ik posed the following problem.
{
\renewcommand{\thetheorem}{28}
\begin{problem}
What is the maximum Wiener index among graphs of order $n$ and diameter $d$?
\end{problem}
}
An asymptotic solution to the problem of Plesn\'ik was obtained by Cambie \cite{cambie}.

\begin{theorem}\label{cam}
There exist positive constants $c_1$ and $c_2$ such that for any $d\geq 3$ the following holds. The maximum Wiener index among all graphs of diameter $d$ and order $n$ is between $d-c_1\frac{d^{3/2}}{\sqrt{n}}$ and $d-c_2\frac{d^{3/2}}{\sqrt{n}}$, i.e. it is of the form $d-\Theta \left ( \frac{d^{3/2}}{\sqrt{n}} \right )$.
\end{theorem}

While the general problem of Plesn\' ik remains open, several significant advances have been made. Wang and Guo \cite{WanGuo} were among the first to determine the trees with the maximum Wiener index for certain specific cases $-$ namely, when the diameter $d$ satisfies $2 \leq d \leq 4$ or $n-3 \leq d \leq n-1$, where $n$ is the order of the tree. Independently, Mukwembi and Vetrík \cite{MukVet} solved the problem for trees with diameter up to 6. 
In addition, Cambie extends their result  for trees in \cite{cambie} by strengthening  the  upper bound given in Theorem \ref{cam}. We also note that trees of order $n$ and diameter $d$ minimizing the Wiener index were characterized in \cite{LiuPan}. For a detailed survey covering these and other selected topics related to the Wiener index, we refer the reader to \cite{tepeh}. \\

DeLaViña and Waller \cite{DelWal} proposed the following conjecture with further  restrictions in problem posed by Plesn\'ik.

\begin{conjecture}
	Let $G$ be a graph with diameter $d>2$ and order $2d+ 1$. Then $W(G)\leq W(C_{2d+1})$, where $C_{2d+1}$ denotes the cycle of length $2d+ 1$.
\end{conjecture}

Motivated by the above conjecture and the results obtained in \cite{skreko},  in this paper we consider the class of trees with diameter $d$ and order $$n\leq d-2+ \left \lfloor \frac{d-1}{2} \right \rfloor$$
and identify those that maximize the Wiener index.  
In \cite{skreko} the authors prove that a graph on $n$ vertices with diameter $d=n-c$ where  $$n\geq \frac 16(7c^3-18c^2+23c-6)$$
has maximum Wiener index if and only if it is a balanced double broom graph. If the above inequality is fulfilled then  $$c\leq \sqrt[3]{\frac 67d+\frac 67}.$$ 
We prove that in the class of trees (with large enough $d$) a double broom graph will maximize the Wiener index if 
$$c\leq 4 \left\lfloor \sqrt{ \frac{d-1}{2}} \right\rfloor-2$$
thereby improving the result given in \cite{skreko} for trees. We also prove that this bound is essentially tight,  more precisely,  if 
$$c\geq 4 \left\lfloor \sqrt{ \frac{d-1}{2}} \right\rfloor+4,$$ 
then a double broom of diameter $d$ and order $d+c$ does not maximize the Wiener index (in the class of trees with diameter $d$ and order $d+c$).

\section{Results}

Let $n$ and $d$ be integers such that $n>d$ and $\mathcal T$ the class of trees of  order $n$ and 
 diameter $d$. 
When we say that a tree $T$ of order $n$ and diameter $d$ has maximum Wiener index, we mean that $W(T')\leq  W(T)$ for every $T'\in \mathcal T$.   

\begin{lemma}\label{sve-se-može}
Let $T$ be a tree of order $n$ and diameter $d$ with maximum Wiener index and $x$ a leaf of $T$. Then $\ecc(x)=d$. 

\end{lemma}
\proof 
Suppose to the contrary, that $\ecc(x)<d$ for a leaf $x$ of $T$. Let $y$ be the nearest vertex to $x$ of degree at least 3 in $T$. 
Let $y'$ and $y''$ be neighbors of $y$ that are not contained in the $xy$-path of $T$, and let $x'$ be the neighbor of $y$ contained in the $xy$-path of $T$ (possibly $x'=x$). If we delete the edge $x'y$ and add the $x'y'$ or $x'y''$ we obtain trees $T_1$ or $T_2$, respectively. 
Since $\ecc(x)<d$ we find that $T_1$ and $T_2$ both have diameter $d$ (note that   $\ecc_{T_i}(x)\leq \ecc_T(x)+1\leq d$). 
If the component of $T-y$ containing $y'$  contains at least as many vertices as the component containing $y''$ then 
$W(T_2)>W(T)$, otherwise $W(T_1)>W(T)$, a contradiction. \qed

Observe that for any leaves $y_1$ and $y_2$ in a tree $T$ of order $n$ and diameter $d$, with maximum Wiener index,  we have $d(y_1,y_2)=d$ or $d(y_1,y_2)$ is even  (this follows from the above lemma).

We call  a vertex $x$ of $T$ a {\em special vertex} if $\deg(x)\geq 3$ and there exist components $T_1$ and $T_2$ of $T-\{x\}$ such that each of them contains exactly one broom vertex of $T$, and for any leaves $y_1\in V(T_1)$ and $y_2\in V(T_2)$ we have $d(y_1,y_2)<d$.   
Hence, by the above observation,  $d(y_1,y_2)$ is even and therefore $d(x,y_1)=d(x,y_2)$.

\begin{definition} \label{miki} 
Let $T$ be tree,  $x$  a special vertex of $T$, and let $T_1$ and $T_2$ be two  
 components of $T-x$ such  that each of them contains exactly one broom vertex of $T$. 
Let  $y_1\in V(T_1)$ and $y_2\in V(T_2)$ be leaves of  $T$ and denote by $y_1'$ and $y_2'$  the broom vertices adjacent to $y_1$ and 
$y_2$, respectively. 
Let $T'$ be the tree obtained from $T$ by deleting $T_2$ and attaching $|V(T_2)|$ leaves to $y_1'$. 
We denote the described operation (of obtaining $T'$ from $T$) by  $T_2\xrightarrow{T} y_1'$, and we write  
 $T'=T_2\xrightarrow{T} y_1'$.  We call the described transformation {\em "relocating a broom"}.
\end{definition}

\begin{figure}[ht!]
	\begin{center}
		\begin{tikzpicture}[scale=.6,style=thick,x=1cm,y=1cm]
			\def\vr{2.5pt} 
			\def\dr{0.5pt}
			\node at (5,-1){$\vdots$};
			\path (5,0) coordinate (x);
			\path (6,1) coordinate (v1);
			\path (8,3) coordinate (v2);
			\path (9,4) coordinate (v3);
			\path (6.8,1.8) coordinate (p1);
			\path (7,2) coordinate (p2);
			\path (7.2,2.2) coordinate (p3);
			\path (3.2,1.8) coordinate (p4);
			\path (3,2) coordinate (p5);
			\path (2.8,2.2) coordinate (p6);
			
			\path (4,1) coordinate (x1);
			\path (2,3) coordinate (x2);
			\path (1,4) coordinate (x3);
			
			\path (-0.5,5.5) coordinate (x4);
			\path (0,5.5) coordinate (x5);
			\path (0.5,5.5) coordinate (x6);
			\node at (1.25,5.5){$\dots$};
			\path (2,5.5) coordinate (x7);
			
			\node at (0.75,6.5){$T_1$};
			\node at (8.75,6.5){$T_2$};
			
			\path (7.5,5.5) coordinate (v4);
			\path (8,5.5) coordinate (v5);
			\path (8.5,5.5) coordinate (v6);
			\node at (9.25,5.5){$\dots$};
			\path (10,5.5) coordinate (v7);
			
			
			
			\draw (x) -- (v1);
			\draw (v2) -- (v3);
			\draw (v3) -- (v4);
			\draw (v3) -- (v5);
			\draw (v3) -- (v6);
			\draw (v3) -- (v7);
			\draw (x) -- (x1);
			\draw (x2) -- (x3);
			\draw (x4) -- (x3);
			\draw (x5) -- (x3);
			\draw (x6) -- (x3);
			\draw (x7) -- (x3);
			
			\draw (x) [fill=white] circle (\vr);
			\draw (v1) [fill=red] circle (\vr);
			\draw (v2) [fill=red] circle (\vr);
			\draw (v3) [fill=red] circle (\vr);
			\draw (v4) [fill=red] circle (\vr);
			\draw (v5) [fill=red] circle (\vr);
			\draw (v6) [fill=red] circle (\vr);
			\draw (v7) [fill=red] circle (\vr);
			\draw (x1) [fill=black] circle (\vr);
			\draw (x2) [fill=black] circle (\vr);
			\draw (x3) [fill=black] circle (\vr);
			\draw (x4) [fill=black] circle (\vr);
			\draw (x5) [fill=black] circle (\vr);
			\draw (x6) [fill=black] circle (\vr);
			\draw (x7) [fill=black] circle (\vr);
			\draw (p1) [fill=black] circle (\dr);
			\draw (p2) [fill=black] circle (\dr);
			\draw (p3) [fill=black] circle (\dr);
			\draw (p4) [fill=black] circle (\dr);
			\draw (p5) [fill=black] circle (\dr);
			\draw (p6) [fill=black] circle (\dr);
			
			\draw[anchor = north] (x) node {$x$};
			\draw[anchor = north] (v3) node {$y_2'$};
			\draw[anchor = south] (v7) node {$y_2$};
			\draw[anchor = north east, inner sep=1pt] (x3) node {$y_1'$};
			\draw[anchor = south] (x7) node {$y_1$};
			
			\node at (13.5,1.75){$\longrightarrow$};
			
			\node at (19.5,-1){$\vdots$};
			\path (19.5,0) coordinate (xa);
			
			\path (19.5,1) coordinate (x1a);
			\node at (19.5,2.1){$\vdots$};
			\path (19.5,3) coordinate (x2a);
			\path (19.5,4) coordinate (x3a);
			
			\path (16.5,5.5) coordinate (x4a);
			\path (17,5.5) coordinate (x5a);
			\path (17.5,5.5) coordinate (x6a);
			\node at (18.25,5.5){$\dots$};
			\path (19,5.5) coordinate (x7a);
			
			
			\path (20,5.5) coordinate (v4a);
			\path (20.5,5.5) coordinate (v5a);
			\path (21,5.5) coordinate (v6a);
			\node at (21.75,5.5){$\dots$};
			\path (22.5,5.5) coordinate (v7a);
			
			\draw[decorate,decoration={brace,amplitude=10pt}]
			(19.8,5.7) -- (22.7,5.7) node[midway, above=10pt] {$|V(T_2)|$};
			
			
			\draw (x3a) -- (v4a);
			\draw (x3a) -- (v5a);
			\draw (x3a) -- (v6a);
			\draw (x3a) -- (v7a);
			\draw (xa) -- (x1a);
			\draw (x2a) -- (x3a);
			\draw (x4a) -- (x3a);
			\draw (x5a) -- (x3a);
			\draw (x6a) -- (x3a);
			\draw (x7a) -- (x3a);
			
			\draw (xa) [fill=white] circle (\vr);
			\draw (v4a) [fill=red] circle (\vr);
			\draw (v5a) [fill=red] circle (\vr);
			\draw (v6a) [fill=red] circle (\vr);
			\draw (v7a) [fill=red] circle (\vr);
			\draw (x1a) [fill=black] circle (\vr);
			\draw (x2a) [fill=black] circle (\vr);
			\draw (x3a) [fill=black] circle (\vr);
			\draw (x4a) [fill=black] circle (\vr);
			\draw (x5a) [fill=black] circle (\vr);
			\draw (x6a) [fill=black] circle (\vr);
			\draw (x7a) [fill=black] circle (\vr);
			
			\draw[anchor = north] (xa) node {$x$};
			\draw[anchor = north west, inner sep=1pt] (x3a) node {$y_1'$};
			\draw[anchor = south] (x7a) node {$y_1$};

		\end{tikzpicture}
	\end{center}
	\caption{The operation $T_2\xrightarrow{T} y_1'$ from Definition \ref{miki}, where black vertices belong to $T_1$ and red vertices belong to $T_2$.}
\end{figure}
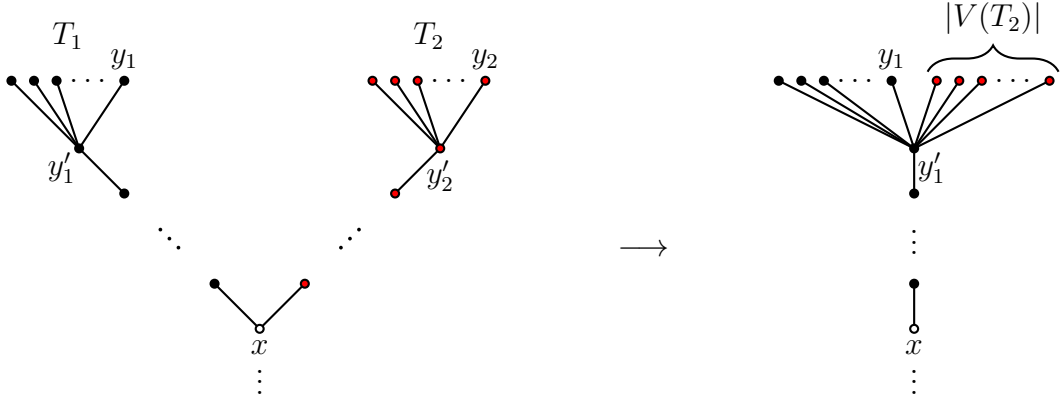
 
In the following lemmas starting with Lemma \ref{optimist} and ending with Corollary \ref{okok} we use 
the notation  as in Definition \ref{miki} for 
$x, T_1, T_2,y_1,y_1',y_2$ and $y_2'$. Additionaly, we denote $p=d(x,y_1)=d(x,y_2)$ and we let $t_i$ be  the number of leaves of $T$ contained in $T_i$ for $i\in [2]$.  

\begin{lemma}\label{optimist}
Let $T$ be a tree of order $n$ and diameter $d$ with maximal Wiener index. 
Let $A=V(T_1)\cup \{x\}$ and $B=V(T)\setminus (A\cup V(T_2))$. If $T'=T_2\xrightarrow{T} y_1'$, 
 then 
\begin{eqnarray} \label{pozitivista}
&&W(T')-W(T) = (p+t_2-1)(p+t_2-2)-\frac{p(p-1)}{6}(3t_2+p-2)-t_2(t_2-1)+\nonumber\\
&& |A\cup B|(p+t_2-1-\frac p2(p-1+2t_2))+ (p+t_2-1)(p-1)(|B|-t_1).
\end{eqnarray}

\end{lemma}


\proof Let $\overline{T_2}=T\backslash T_2$ and $\overline{T_2'}=T'\backslash T_2'$, where $T_2'$ denotes the subtree in $T'$ induced by vertices of $V(T_2)$. Then
$$W(T)=W(T_2)+W(\overline{T_2})+d(T_2, \overline{T_2}),$$
\noindent where $d(T_2, \overline{T_2})$ represents the sum of distances between vertices of $T_2$ and $\overline{T_2}$.
For tree $T'$ it holds
$$W(T')=W(T_2')+W(\overline{T_2'})+d(T_2', \overline{T_2'}),$$
\noindent where $d(T_2',\overline{T_2'})$ is the sum of distances between vertices of $T_2'$ and $\overline{T_2'}$.

\noindent It can easily be concluded that
\begin{eqnarray*}
& & W(T_2)={p \choose 3}+{t_2\choose 2}\cdot 2+t_2 \sum_{i=1}^{p-1} i=\frac{p(p-1)}{6}(3t_2+p-2)+t_2(t_2-1),\\
& & W(T_2')={p+t_2-1\choose 2}\cdot 2=(p+t_2-1)(p+t_2-2),\\
& & W(\overline{T_2})=W(\overline{T_2'}).
\end{eqnarray*}

\noindent In the next calculations we use the fact that $|\overline{T_2}|=|\overline{T_2'}|=|A\cup B|$ and $|T_2|=|T_2'|=p+t_2-1$. \\
\noindent By the definition of $d(T_2,\overline{T_2})$ we have
\begin{eqnarray*}
d(T_2,\overline{T_2})&=&\sum\limits_{k=1}^{p-1}\left(\sum\limits_{u\in A}(k+d(u,x))+\sum\limits_{u\in B}(k+d(u,x))\right)\\
&+&t_2\left(\sum\limits_{u\in A}(p+d(u,x))+\sum\limits_{u\in B}(p+d(u,x))\right)\\
&=&\sum_{k=1}^{p-1} \sum_{u\in A\cup B} k + \sum_{k=1}^{p-1} \sum_{u\in A\cup B} d(u,x) +t_2\sum_{u\in A\cup B} p + t_2\sum_{u\in A\cup B} d(u,x)\\
&=& \left(\frac{p(p-1)}{2}+t_2p\right)\cdot|A\cup B| + (p+t_2-1) \sum_{u\in \overline{T_2}} d(u,x) \\
&=&\frac{p}{2}(p-1+2t_2)\cdot |A\cup B|+(p+t_2-1)\sum\limits_{u\in\overline{T_2}}d(u,x).
\end{eqnarray*}
Let $y_1'$ be the broom vertex adjacent to leaves of $T_1$. Then
\begin{eqnarray*}
d(T_2',\overline{T_2'})&=&(p+t_2-1)\sum\limits_{u\in\overline{T_2'}}(1+d(u,y_1'))=\\
&=&(p+t_2-1)|\overline{T_2'}|+(p+t_2-1)\sum\limits_{u\in\overline{T_2'}}d(u,y_1')\\
&=&(p+t_2-1)|A\cup B|+(p+t_2-1)\sum_{u\in\overline{T_2'}}d(u,y_1').\end{eqnarray*}

\noindent Now we have
\begin{eqnarray*}
	W(T')-W(T)&=&W(T_2')-W(T_2)+d(T_2', \overline{T_2'})-d(T_2,\overline{T_2})\\
	&=&(p+t_2-1)(p+t_2-2)-\frac{p(p-1)}{6}(3t_2+p-2)-t_2(t_2-1)\\
	&+&|A\cup B|\left(p+t_2-1-\frac{p}{2}(p-1+2t_2) \right)\\
	&+&(p+t_2-1)\left (\sum_{u\in A\cup B}(d(u,y_1')-d(u,x)) \right )\\
	&=&(p+t_2-1)(p+t_2-2)-\frac{p(p-1)}{6}(3t_2+p-2)-t_2(t_2-1)\\
	&+&|A\cup B|\left(p+t_2-1-\frac{p}{2}(p-1+2t_2)\right)+ \\
	&+&(p+t_2-1) \left ( \underbrace{\sum_{u\in A}(d(u,y_1')-d(u,x))}_{-(p-1)t_1} + \underbrace{\sum_{u\in B}(d(u,y_1')-d(u,x))}_{(p-1)|B|} \right ) \\ 
	&=&(p+t_2-1)(p+t_2-2)-\frac{p(p-1)}{6}(3t_2+p-2)-t_2(t_2-1)\\
	&+&|A\cup B|\left(p+t_2-1-\frac{p}{2}(p-1+2t_2)\right)+(p+t_2-1)(p-1)(|B|-t_1).
\end{eqnarray*}
\qed

Using  the notation of  Definition \ref{miki} we have the following corollaries.

\begin{corollary} \label{nepobjedljivi}
Let $T$ be a tree of diameter $d$ on $n$ vertices with maximum Wiener index, and let $T'=T_2\xrightarrow{T} y_1'$. 
Then 
\begin{eqnarray*}
&&W(T')-W(T) = -\frac 16(p-1)(12(t_1-1)(p+t_2-1)+p(-3n-5+12t_2+10p)).
\end{eqnarray*}
\end{corollary}
\proof Use $|A\cup B|=n-(p+t_2-1)$ and $|B|=n-(p+t_2-1)-(p+t_1)=n-2p-t_1-t_2+1$ in \eqref{pozitivista} 
to obtain the claimed result. \qed

The following corollary is a straightforward consequence of Corollary \ref{nepobjedljivi}.

\begin{corollary}\label{optimist8}
Let $T$ be a tree of diameter $d$ on $n$ vertices with maximum Wiener index, and let  $T'=T_2\xrightarrow{T} y_1'$. Then 
$W(T)\geq W(T')$  if and only if
\begin{eqnarray*}
&& t_1-1 \geq \frac{p(3n+5- 12t_2-10p)}{12(p+t_2-1)}.
\end{eqnarray*}
\end{corollary}

\begin{corollary}\label{optimist-forever}
Let $T$ be a tree of diameter $d$ on $n$ vertices with maximum Wiener index, and let $T'=T_2\xrightarrow{T} y_1'$.  If  $t_1=t_2$, then 
$W(T)\geq W(T')$  if and only if
\begin{eqnarray} \label{veliki-optimist}
 t_1 \geq  \sqrt{\frac{p(3n+2p-7)}{12}}-p+1.
\end{eqnarray}
\end{corollary}

\proof Use  $t_1=t_2$ in Corollary \ref{optimist8}. \qed

\begin{corollary}\label{najjači-optimist}
Let $T$ be a tree of diameter $d$ on $n$ vertices with maximum Wiener index, and let $T'=T_2\xrightarrow{T} y_1'$.  If  $t_1=t_2+1$, then 
$W(T)\geq W(T')$  if and only if
\begin{eqnarray}\label{največji-optimist}
&& t_1 \geq  \sqrt{\frac{p(3n+2p-7)+3}{12}}-p+\frac 32.
\end{eqnarray}

\end{corollary}

\proof Use  $t_1=t_2+1$ in Corollary \ref{optimist8}. \qed

\begin{lemma}\label {nono}
If $n\geq 1636$, then $f(p)=\sqrt{\frac{p(3n+2p-7)}{12}}-p+1$ and $g(p)=\sqrt{\frac{p(3n+2p-7)+3}{12}}-p+\frac 32$ are   increasing functions of $p$ 
for $2\leq p \leq 4\sqrt{\frac{n-1}{2}}-3$.
\end{lemma}

\proof 
The derivative of $f$ is 
$$f'(p)=\frac{1}{2\sqrt{12}}\frac{3n+4p-7}{\sqrt{p(3n+2p-7)}}-1> 
\frac{1}{2\sqrt{12}} \sqrt{\frac{3n+4p-7}{p}}-1.$$
Hence $f'(p)> 0$, if
$$\sqrt{\frac{3n+4p-7}{p}}\geq 2\sqrt{12}$$
which simplifies to $3n-7\geq 44p$.
Since  $p \leq 4\sqrt{\frac{n-1}{2}}-3$ we have $f'(p)\geq 0$ whenever 
$$3n-7\geq 176\sqrt{\frac{n-1}{2}}-132.$$
Solving this quadratic inequality gives $n\geq 1636$. Similar proof works for the function $g$. 
\qed

The following corollary is straightforward (we use $p=2$ in Corollary \ref{optimist-forever}
 and \ref{najjači-optimist} and Lemma \ref{nono})

\begin{corollary}\label{okok}
Let $T$ be a tree of diameter $d$ on $n$ vertices with maximum Wiener index, and let $T'=T_2\xrightarrow{T} y_1'$.  If   
$n\geq 1636, p \leq 4\sqrt{\frac{n-1}{2}}-3$ and 
$W(T)\geq W(T')$  then
 \begin{eqnarray*}
&& t_1 \geq  \left \lceil\sqrt{\frac{n-1}{2}}\right\rceil-1 ~{\rm if}~ t_1=t_2,~{\rm and }\\
&& t_1 \geq  \left \lceil\sqrt{\frac{n-\frac 12}{2}}-\frac 12\right\rceil ~{\rm if} ~t_1=t_2+1.
\end{eqnarray*} 
\end{corollary}

The following lemma gives the change of Wiener index when relocating a leaf of the tree.

\begin{lemma}\label{optimist-sve-može}
Let  $T$ be a tree and $x,y$ leaves of $T$. Let $P$ be the $xy$-path in $T$ and $x',y'$ vertices of $T$ adjacent to $x,y$, respectively. 
 Let $V(P)=[r]_0$ where  $x=0$ and $y=r$.  
Let $T_i$ be the union of components of $T-P$ adjacent to $i$ in $T$ for $i\in [r-1]$ (if there is no such component then let  $T_i$ be the empty graph).
If $T'$ is the tree obtained from $T$ by deleting $x$ and adding a leaf to $y'$, then  
 \begin{equation}\label{vječiti-optimist}
W(T')-W(T)=\sum_{i=1}^{r-1}(r-2i)|V(T_i)|-(r-2).
\end{equation}
\end{lemma}

\proof Let $z\in V(T')\setminus V(T)$ be the leaf of $T'$ adjacent to $y'$. 
Then we have $W(T')=W(T-x)+d(z,T-x)$ and $W(T)=W(T-x)+d(x,T-x)$. 
Now we have 
\begin{eqnarray} \label{još-ima-nade-za-nas}
  &&d(x,T-x)=d(x,P-x)+\sum_{i=1}^{r-1}d(x,T_i)=d(x,P-\{x,y\})+d(x,y)+\\
\nonumber &&\sum_{i=1}^{r-1}d(x,i)|V(T_i)|+\sum_{i=1}^{r-1}d(i,T_i)=d(x,P-\{x,y\})+r+\sum_{i=1}^{r-1}i|V(T_i)|+\sum_{i=1}^{r-1}d(i,T_i)
\end{eqnarray}

and similarly
\begin{eqnarray} \label{svakog-dana-sve-više-napredujem}
  &&d(z,T-x)=d(z,P-x)+\sum_{i=1}^{r-1}d(z,T_i)=d(z,P-\{x,y\})+d(z,y)+\\
\nonumber &&\sum_{i=1}^{r-1}d(z,i)|V(T_i)|+\sum_{i=1}^{r-1}d(i,T_i)=d(z,P-\{x,y\})+2+\sum_{i=1}^{r-1}(r-i)|V(T_i)|+\sum_{i=1}^{r-1}d(i,T_i).
\end{eqnarray}
Since $W(T')-W(T)=d(z,T-x)-d(x,T-x)$ we obtain \eqref{vječiti-optimist} by subtracting  \eqref{još-ima-nade-za-nas} from \eqref{svakog-dana-sve-više-napredujem}. 
\qed

The transformation in which $T'$ is obtained from $T$, as described in Lemma \ref{optimist-sve-može}, is called {\em relocating a leaf}.  
It is straightforfard to prove (by an application of  Lemma \ref{optimist-sve-može}) 
that in the class of double brooms with a fixed order and diameter the balanced double broom maximizes the Wiener index.  

\begin{corollary}\label{mi-smo-nevjerovatni}
Let  $T$ be a tree on $n$ vertices with maximum Wiener index and let  $x$ be a special vertex of $T$. 
Let $T_1$ and $T_2$ be    
 components of $T-x$ such  that each of them contains exactly one broom vertex of $T$. 
For $i\in [2]$ let  $t_i$ be  the number of leaves of $T$ contained in $T_i$. 
Then $|t_1-t_2|\leq 1$. 
\end{corollary}
\proof 
Let $y_1\in V(T_1)$ and $y_2\in V(T_2)$ be two leaves of $T$. Let   $p=d_T(x,y_1)=d_T(x,y_2)$.  If $T'$ and $T''$ are trees obtained from $T$ by moving 
$y_1$ and $y_2$ to $T_2$ and $T_1$, respectively, then by \eqref{vječiti-optimist}  we have 
\begin{eqnarray*}
W(T')-W(T)&=&
(2p-2)(t_1-t_2-1) ~~{\rm and }\\
W(T'')-W(T)&=&
(2p-2)(t_2-t_1-1).
\end{eqnarray*}
It follows that $t_1-t_2-1\leq 0$ and $t_2-t_1-1\leq 0$,   hence   $|t_1-t_2|\leq 1$. 
\qed

\begin{lemma}\label{kdaj-bo-ta-clanek-fertig?}
Let  $T$ be a tree of diameter $d$ on $n$ vertices with maximum Wiener index. Let $P$ be a path of length $d$ in $T$, and 
$x$ a special vertex of $P$. Assume that $T_1$ and $T_2$ are  components of $T-x$ having   exactly one broom vertex. Then 
the number of vertices in $T_1\cup T_2$ that are not in $P$ is at least 
 $$\left \lceil2\sqrt{\frac{n-1}{2}}\right\rceil-2.$$
\end{lemma}

\proof
Let $b_1$ and $b_2$ be the broom vertices in $T_1$ and $T_2$, and let $t_1$ and $t_2$ be the number of leaves attached to 
$b_1$ and $b_2$, respectively. 
By Corollaries \ref{optimist-forever} and \ref{najjači-optimist}
$$t_1 +t_2\geq \left \lceil2\sqrt{\frac{p(3n+2p-7)}{12}}\right\rceil-2p+2,$$
where $2p-2=d_T(b_1,b_2)$. Note that $p\leq d/2$ 
and that at most one leaf of $T_1\cup T_2$ is contained in $P$. Let $P_1$ be the path from $x$ to $b_1$ and $P_2$ the path from $x$ to $b_2$ and note that  
either $P_1-x$ or $P_2-x$ is disjoint with $P$. Since $|P_1|=|P_2|=p$ it follows that 
the number of vertices in $T_1\cup T_2$ that are not in $P$ is at least 
$$\left \lceil2\sqrt{\frac{p(3n+2p-7)}{12}}\right\rceil-2p+2-1+(p-1)=\left \lceil2\sqrt{\frac{p(3n+2p-7)}{12}}-p\right\rceil.$$
The above function is  increasing for $2\leq p\leq d/2\leq n/2$ (one may verify this by showing that the derivative of the function
$f(p)=2\sqrt{\frac{p(3n+2p-7)}{12}}-p$ is nonegative for  $2\leq p\leq n/2$). 
It follows that the above expression is minimum when $p=2$ and it is 
 $$\left \lceil2\sqrt{\frac{n-1}{2}}\right\rceil-2.$$\qed

\begin{theorem}
Let $T$ be a tree of diameter $d$ on $n\leq d-2+4 \left\lfloor \sqrt{ \frac{d-1}{2}} \right\rfloor$ vertices. 
If $n\geq 1636$ then $W(T)$ is maximum if and only if  $T$ is a  balanced double broom. 
\end{theorem}

\proof
Suppose that $T$ is a tree with diameter $d$ on $n\leq d-2+4 \left\lfloor \sqrt{ \frac{d-1}{2}} \right\rfloor$ vertices with maximum Wiener index, where $n\geq 1636$.
We will prove that if $T$ is not a double broom graph 
then there is a transformation (either relocating a broom or relocating a leaf) on $T$, that doesn't change the diameter and the order of $T$, which increases the Wiener index, contradicting the assumption that $T$ has maximum Wiener index. Hence 
double brooms are the only candidates for the maximum Wiener index of a tree (subject to fixed $d$ and $n$). Moreover, by  Lemma \ref{optimist-sve-može}, 
a double broom that maximizes the Wiener index is balanced --- this implies the theorem.

 Let 
$P$ be a path of length $d$ in $T$ and assume that $V(P)=[d]_0$. Note that if $T$ has no special vertex, then it is a double broom, and we are done. 
So assume that $T$ has at least one special vertex, and we may assume that a special vertex of $T$ is contained in $P$ (since we can choose $P$ so that this is the case).  Since  $n\leq d-2+4 \left\lfloor \sqrt{ \frac{d-1}{2}} \right\rfloor$  
we find that there are at most  $4 \left\lfloor \sqrt{ \frac{d-1}{2}} \right\rfloor-3$ vertices in $V(T)\setminus V(P)$.

We claim that there is at most one special vertex of $T$. Assume the contrary, that there are special vertices $x$ and $y$ of $T$. Let 
$T_1, T_2$ and $R_1,R_2$ be components of $T-x$ and $T-y$, respectively, having exactly one broom vertex. Then 
by Lemma \ref{kdaj-bo-ta-clanek-fertig?}  the number of vertices of 
$T_1\cup T_2\cup R_1\cup R_2$ that are not contained in $P$ is at least $$\displaystyle2\left \lceil2\sqrt{\frac{n-1}{2}}\right\rceil-4$$ (note that $T_1\cup T_2$ is disjoint with $R_1\cup R_2$) 
and therefore    
  $$|V(T)\setminus V(P)|\geq 2\left \lceil2\sqrt{\frac{n-1}{2}}\right\rceil-4>4 \left\lfloor \sqrt{ \frac{d-1}{2}} \right\rfloor-3$$
(note that if $a>b$ then $\lceil{2a}\rceil>2\lfloor b\rfloor$), a contradiction. 
This proves the claim.

Assume therefore that $T$ has  exactly one special vertex, say vertex $p\leq \lceil d/2\rceil-1$. 
Let $T_p$ be a component of $T-p$ disjoint with $P$ having exactly one broom vertex. 
Let $t_1, t_{d-1}$ and $t_{p}$ be the number of leaves 
adjacent to vertices $1, d-1$ and the number of leaves in $T_p$, respectively. 
Since $T$ has exactly one special vertex we find that every vertex $u\geq  \lceil d/2\rceil, u\neq d-1$ on $P$ is of degree 2 in $T$ 
(for otherwise $T$ would surely have a special vertex distinct from $p$).  

Let $T'$ be the tree obtained from $T$ by relocating a leaf adjacent to vertex $1$ 
and attaching it to vertex $d-1$. For any vertex $u\in V(P)$ we denote the union  of components of $T-u$ that are disjoint with $P$ by $T_u$, and we note that 
$T_u=\emptyset$ if $u\geq  \lceil d/2\rceil, u\neq d-1$. Thus, 
by Lemma \ref{optimist-sve-može}, we have  

\begin{equation}\label{optimist2}
W(T')-W(T)\leq t_1(d-2)+(t_p+p-1)(d-2p)-t_{d-1}(d-2)-(d-2)\leq 0
\end{equation}
where the last inequality follows from maximality of $T$.  
The last inequality in \eqref{optimist2} simplifies to 
\begin{equation}\label{zzz} 
  (t_1+t_p)(d-2)+(p-1)(d-2-2(t_p+p-1))\leq (t_{d-1}+1)(d-2).
\end{equation}
Since $t_p+p-1\leq 4 \sqrt{ \frac{d-1}{2}} -3$  we find that $d-2-2(t_p+p-1)>0$ for $d\geq 22$ (which we know to be the case, since $n\geq 1636$),  
and therefore, by \eqref{zzz} we have 
  $t_{d-1}\geq t_1+t_p$.

By Corollary \ref{mi-smo-nevjerovatni} we have  $|t_p-t_1|\leq 1$. 
Since $n\geq 1636$ we have,
by Corollary \ref{okok}, 
 \begin{eqnarray*}
&& t_1 \geq  \left \lceil\sqrt{\frac{n-1}{2}}\right\rceil-1 ~{\rm if}~ t_1=t_p,~{\rm and }\\
&& t_1 \geq  \left \lceil\sqrt{\frac{n-\frac 12}{2}}-\frac 12\right\rceil ~{\rm if} ~t_1=t_p+1.
\end{eqnarray*} 
We now distinguish two cases: \\ 

{\em Case 1:} Suppose that $t_1=t_p$. 
Then 

\begin{equation*}\label{optimist4}
|V(T)\setminus V(P)|\geq t_1+t_p+t_{d-1}-1\geq 4\left\lceil \sqrt{\frac{n-1}{2}}\right\rceil -5> 4 \left\lfloor \sqrt{ \frac{d-1}{2}} \right\rfloor-3,
\end{equation*}
which is a contradiction.

{\em Case 2:} Suppose that $t_1=t_p+1$. 
Then 

\begin{equation}\label{dajže}
 |V(T)\setminus V(P)|\geq t_1+t_p+t_{d-1}-1\geq 2\left (\left\lceil \sqrt{\frac{n-\frac 12}{2}}-\frac 12\right \rceil+\left\lceil \sqrt{\frac{n-\frac 12}{2}}-\frac 32\right\rceil\right) -1
\end{equation}
Now assume first that $n\leq  d-3+4 \left\lfloor \sqrt{ \frac{d-1}{2}}\right\rfloor$. Then there are at most  $4 \left\lfloor \sqrt{ \frac{d-1}{2}} \right\rfloor-4$ vertices in $V(T)\setminus V(P)$. 
Note that for any integer $a$ we have  $\lceil a-\frac 12\rceil+\lceil a-\frac 32\rceil\geq 2\lceil a\rceil-3$, and therefore by \eqref{dajže} we have

\begin{equation*}\label{dajže1}
 |V(T)\setminus V(P)|\geq 4\left\lceil \sqrt{\frac{n-\frac 12}{2}}\right \rceil -7>4 \left\lfloor \sqrt{ \frac{d-1}{2}} \right\rfloor-4,
\end{equation*}
a contradiction.
Suppose therefore that  $n= d-2+4 \left\lfloor \sqrt{ \frac{d-1}{2}}\right\rfloor$.  It is easy to see that then 
$$\left\lceil \sqrt{\frac{n-\frac 12}{2}}-\frac 32\right\rceil\geq  \left\lfloor \sqrt{ \frac{d-1}{2}} \right\rfloor$$ 
(one can prove this by proving that $ \sqrt{\frac{n-\frac 12}{2}}-\frac 32>   \sqrt{ \frac{d-1}{2}}-1$ for $d\geq 9$), and therefore by \eqref{dajže} we have 
$$  |V(T)\setminus V(P)|\geq  4 \left\lfloor \sqrt{ \frac{d-1}{2}} \right\rfloor+1,$$
a contradiction.  \qed

\begin{lemma} \label{najvecjidoublebroom}
	Let $n,a$ and $b$ be positive integers and $g=a+b$. Then the Wiener index of $B(n,a,b)$ is maximum if $a=b=g/2$ for even $g$ and 
	$a=(g - 1)/2, b=(g + 1)/2$ for odd $G$.
\end{lemma}

\proof First, assume that $g$ is even, and let $B=B\left(n, \frac{g}{2}, \frac{g}{2}\right)$. Let $B'=B'(n,\frac{g}{2}-1, \frac{g}{2}+1)$ be the double broom obtained from $B$ by deleting an arbitrary leaf neighbor of $x$
and adding a leaf to $y$. By Lemma \ref{optimist-sve-može} it holds that
$$W(B')-W(B)=(r-2)\frac{g}{2}+(2-r)\frac{g}{2}-r+2=2-r<0.$$

Now assume that $g$ is odd, and let $B=B\left(n, \frac{g-1}{2}, \frac{g+1}{2}\right)$. Let $B'=B'(n,\frac{g-1}{2}-1, \frac{g+1}{2}+1)$, again obtained from $B$ by deleting an arbitrary leaf neighbor of $x$
and adding a leaf to $y$. Then, by Lemma \ref{optimist-sve-može}, we get
$$W(B')-W(B)=(r-2)\frac{g-1}{2}+(2-r)\frac{g+1}{2}-r+2=2(2-r)<0.$$

This completes the proof. \qed

\begin{proposition}
If $n\geq  d+4+4 \left\lfloor \sqrt{ \frac{d-1}{2}} \right\rfloor$ 
there is a triple broom $B(n,a,b,c)$ of diameter $d$ such that for any double broom $B(n,q,r)$ of diameter $d$ we have 
$W(B(n,a,b,c))>W(B(n,q,r)).$
\end{proposition}

\proof
We prove the proposition for $n= d+4+4 \left\lfloor \sqrt{ \frac{d-1}{2}} \right\rfloor$; the proof for 
$n> d+4+4 \left\lfloor \sqrt{ \frac{d-1}{2}} \right\rfloor$ is similar. 
Let $T$ be a triple broom denoted by $B(n,a,b,c)$, where $n=d+g$ and 
$$g= 4+4 \left\lfloor \sqrt{ \frac{d-1}{2}} \right\rfloor.$$ 
Let $T'=T_2\xrightarrow{T} x=B(n,q,r)$ be a corresponding double broom. 
This means that $B(n,a,b,c)$ has $g$ leaves, while $B(n,q,r)$ has $g+1$ leaves. By Lemma \ref{najvecjidoublebroom}, it follows that among all possible double brooms on 
$n$ vertices, the maximum Wiener index is attained by $B \left ( n,\frac{g}{2},\frac{g+2}{2} \right )$.
Using Corollary \ref{optimist-forever}, we will show that 
$$W\left ( B \left ( n,\frac{g}{4},\frac{g}{2},\frac{g}{4} \right ) \right ) >W \left ( B \left ( n,\frac{g}{2},\frac{g+2}{2} \right ) \right ).$$ 
Since $t_1=t_2$ and $p=2$, applying inequality (\ref{veliki-optimist}) yields $$\frac g4> \sqrt{\frac{(d+g-1)}{2}}-1,$$ which holds for all $d\geq 3$ provided that $g>\sqrt{8d-24}.$ 
Based on the value of $g$, we observe that $$g\geq 4+4 \left\lfloor \sqrt{ \frac{d-1}{2}} \right\rfloor>4+4\left ( \sqrt{\frac{d-1}{2}}-1 \right )=4 \sqrt{\frac{d-1}{2}}.$$ 
Thus, it suffices to show that $$4 \sqrt{\frac{d-1}{2}}>\sqrt{8d-24}.$$ This inequality is clearly satisfied for all $d\geq 3$, completing the proof. \qed

\section{Concluding remarks}

We have proved that   double broom maximizes the Wiener index in the class of trees with diameter $d$ and order 
$n\leq d-2+4 \left\lfloor \sqrt{ \frac{d-1}{2}} \right\rfloor$.  We think that this remains true in the class of all graphs (if the assumptions on $d$ and $n$ remain the same). Hence we pose  the following question.

\begin{question}
Let $n\leq d-2+4 \left\lfloor \sqrt{ \frac{d-1}{2}} \right\rfloor$. Does the double broom with diameter $d$ and order $n$   maximize Wiener index in the class of all graphs with diameter $d$ and order $n$?
\end{question}

By \cite{cambie3} the above question has a positive answer. 
The conjecture given in \cite{DelWal}, which asserts that the cycle of diameter $d$ and order $2d+1$ maximizes the Wiener index (in the class of graphs with diameter $d$ and order 
$2d+1$),  suggests that when $n$ is sufficiently large relative to $d$, the graph that maximizes the Wiener index is not a tree. 
However, when $n$ is sufficiently small relative to $d$, the results obtained in 
\cite{skreko} imply that a tree (a double broom) will maximize the Wiener index. This brings us to the following problem.

\begin{problem}
Determine the function $f=f(d)$, such that the maximum Wiener index in the class of graphs of order at most $n$ and diameter $d$ is attained by a 
tree if and only if $n\leq f(d)$. 
\end{problem}

\section*{Acknowledgement}
This work  was carried out within the Slovenia–Serbia joint project (for the period 1.7.2023 - 31.12.2025.) entitled “Topological indices of graphs and digraphs”, funded by Slovenian research agency
ARIS, project no. BI-RS/23-25-023, and the Serbian Ministry of Science, Technological Development and Innovation, project no. 337-00-110/2023-05/23. 

The author Bojana Borovi\' canin also gratefully acknowledges additional support from the Ministry of Science, Technological Development and Innovation of the Republic of Serbia, agreement no. 451-03-137/2025-03/ 200122. 

The authors Dragana Bo\v zovi\' c and Simon \v Spacapan were partially supported by the Slovenian Research Agency ARIS: Dragana Bo\v zovi\' c by the research program Telematics P2-0065, and Simon \v Spacapan by project N1-0218 and research program P1-0297. 
The authors thank Stijn Cambie for several helpful comments on our manuscript.

\end{document}